\documentclass{amsart}
\usepackage{amsmath,amssymb,amsthm}
\usepackage[mathscr]{euscript}

\newtheorem*{remark}{Remark}

\newtheorem{thm}{Theorem}

\newtheorem*{goal}{Goal}

\title{Recovering functions from the Paley-Wiener amalgam space $(PW,l^1)$}
\author{Jeff Ledford}
\date{}

\begin{document}

\newcommand{\sinc}{\text{sinc}}
\newcommand{\N}{\mathbb{N}}
\newcommand{\R}{\mathbb{R}}
\newcommand{\Z}{\mathbb{Z}}
\newcommand{\T}{[-\pi,\pi]}
\newcommand{\BV}{\emph{B.V.}}


\begin{abstract}
In this paper we show that functions from the Paley-Wiener amalgam space $(PW,l^1)=\{f\in L^2(\R): \sum\|\hat{f}(\xi+2\pi m) \|_{L^2(\T)} < \infty\}$ enjoy similar recovery properties as the classical Paley-Wiener space.  Specifically, if $\{\phi_\alpha(x): \alpha\in A\}$ is a regular family of interpolators and $\{x_n: n\in\Z\}$ is a complete interpolating sequence for $L^2([-\pi,\pi])$, then the family $\{e^{2\pi i m x}\phi_\alpha(x-x_n):m,n\in\Z, \alpha\in A\}$ may be used to recover $f\in(PW,l^1)$.
\end{abstract}

\maketitle

\section{Introduction}
In this short note, we provide an example of the utility of amalgam spaces, in this case the Paley-Wiener amalgam space
\[
(PW,l^1):=\left\{f\in L^2(\R) : \sum_{m\in\Z}\| \hat{f}(\xi+2\pi m) \|_{L^2(\T)} <\infty \right\}.
\]
Amalgam norms provide both local and global information about a particular function.  One particularly nice feature is that it is often enough to prove a result about the local behavior and use the nature of the mixed norm to extend this to the larger space.  We provide an example of this phenomenon relating to recovery of Paley-Wiener functions using regular family of interpolators.  Most of the heavy lifting is contained in \cite{me}, we only state the salient results found there.  More information on amalgam spaces may be found in \cite{Grochenig}.
The remainder of this note is organized into three sections, the first of which collects definitions and results from \cite{me}, while the second contains the main result and details of the proof. The final section is devoted to examples and conclusions.

\section{Definitions and Basic Facts}
We begin with the convention for the Fourier transform.
For a function $g \in L^1(\R)$, we define the \emph{Fourier transform} of $g$, denoted $\hat{g}$, to be
\begin{equation}\label{F def}
\hat{g}(\xi)=(2\pi)^{-1/2}\int_{\R}g(x)e^{-ix\xi}dx
\end{equation}
We recall the definition of an interpolator, and a family of regualar interpolators.  More information may be found in \cite{me}.
We call $\phi:\mathbb{R}\to\mathbb{R}$ an \emph{interpolator} for $PW_\pi$ if it satisfies satisfies the following conditions:
\begin{enumerate}
\item[(A1)] $\phi(x)\in L^1(\mathbb{R})\cap C(\mathbb{R})$ and $\hat\phi(\xi)\in L^1(\mathbb{R})$.
\item[(A2)] $\hat{\phi}(\xi)\geq 0$ and $\hat{\phi}(\xi)\geq\delta>0$ on $[-\pi,\pi]$.
\item[(A3)] Let $\displaystyle M_j=\sup_{|\xi|\leq\pi}\hat{\phi}(\xi+2\pi j)$, then $\{M_j\}\in l^1(\mathbb{Z})$.
\end{enumerate}
We consider the one parameter family of interpolators $\{\phi_\alpha(x):\alpha\in A\}$, where $\alpha\in A\subset(0,\infty)$.  We will call the family \emph{regular} if it satisfies the following hypotheses.
\begin{itemize}
\item[(H1)] $\phi_{\alpha}(x)$ is an interpolator for $PW_\pi$ for each $\alpha\in A$.
\item[(H2)] $\displaystyle\sum_{j\neq0}M_j(\alpha)\leq C m_\alpha$, where $M_j(\alpha)$ is as in (A3), $\displaystyle m_\alpha = \inf_{|\xi|\leq\pi}\hat{\phi}_{\alpha}(\xi)$, and C is independent of $\alpha$.
\item[(H3)] $\text{For almost every } |\xi|\leq\pi; \displaystyle\lim_{\alpha\to\infty}\dfrac{m_\alpha}{\hat{\phi}_\alpha(\xi)}=0.$
\end{itemize}
These definitions allow us to prove that a certain interpolant $I_\alpha f$ converges to a given Paley-Wiener function $f$.
Specifically, if $\{x_n:n\in\Z\}$ is a complete interpolating sequence for $L^2(\T)$ and we define
\[
I_\alpha f (x) = \sum_{n\in\Z}a_n\phi_\alpha(x-x_n)
\]
where $\{a_n:n\in\Z\}\in l^2$ is chosen so that $I_\alpha f(x_n)=f(x_n)$ for all $n\in\N$, then
\begin{equation}\label{ineq1}
\| f - I_\alpha f  \|_{L^2(\R)} \leq C \| \dfrac{m_\alpha}{\hat\phi_{\alpha}(\xi)} \hat{f}(\xi) \|_{L^2(\T)}
\end{equation}
Here $C$ is a constant independent of both $f$ and $\alpha$.  This inequality may be found in the proof of theorem 1 in \cite{me}.  This formula will be used in the calculations that follow in the next section.

\section{Main Results}
We set ourselves to the following task.
\begin{goal}
Given a function $f\in (PW, l^1)$, we can find a collection of approximands $I_\alpha g_m(x)$ which tend to $f$ uniformly on $\R$.
\end{goal}

\begin{remark}
\emph{While the wording is vague, the idea is to illustrate how to use the results in \cite{me} to prove similar results for the the amalgam space $(PW,l^1)$.  The nature of the approximands will be fleshed out in due course.}
\end{remark}
We begin our procedure by constructing auxiliary functions $g_m \in PW_\pi$ related to $f\in L^2(\R)$.  For $m\in\Z$, let $\hat{f}_m(\xi)= f(\xi)\chi_{\T}(\xi-2\pi m)$.  Now define $g_m$ by its Fourier transform
\[
\hat{g}_m(\xi)=\hat{f}_m(\xi+2 \pi m).
\]
Using the results of \cite{me}, for each $m\in\Z$ there exists a sequence $\{a_{m,n}:n\in\Z\}\in l^2$ such that
\begin{equation}\label{Ig}
I_\alpha g_m(x) = \sum_{n\in\Z}a_{m,n}\phi_\alpha(x-x_n)
\end{equation}
satisfies $I_\alpha g_m(x_n)=g_m(x_n)$ for all $n\in\Z$.  To simplify notation, we write
\[
J_\alpha f (x) = \sum_{m\in\Z}e^{2\pi i m x}I_\alpha g_m(x)
\]
Now we estimate the norm using the same technique that is used in the proof of theorem 1 in \cite{me}.
\begin{align}
&\| f - J_\alpha f (x) \|_{(PW,l^1)} \nonumber   \\
=& \sum_{j\in\Z} \| \hat{f}(\xi+2\pi j) - \sum_{m\in\Z}\widehat{I_\alpha g_m}(\xi-2\pi m + 2\pi j)   \|_{L^2(\T)} \label{eq2} \\
=& \sum_{j\in\Z} \| \hat{g}_j(\xi) - \sum_{m\in\Z}\widehat{I_\alpha g_m}(\xi-2\pi m + 2\pi j)   \|_{L^2(\T)}  \nonumber \\
\leq & C\left(1+\sum_{l\neq 0}M_l(\alpha)m_{\alpha}^{-1}   \right) \sum_{j\in\Z}\| m_\alpha[\hat\phi_\alpha(\xi)]^{-1} \hat{g}_j(\xi)  \|_{L^2(\T)} \nonumber  \\
\leq & C' \sum_{j\in\Z}\| m_\alpha[\hat\phi_\alpha(\xi)]^{-1} \hat{f}(\xi+2\pi j)  \|_{L^2(\T)} \label{bnd1}\\
\leq & C'\|f  \|_{(PW,l^1)} \label{bnd2}
\end{align}
Inequality \eqref{bnd2} allows us to apply the dominated convergence theorem, since all of the terms in \eqref{bnd1} tend to 0 by $(H3)$, we have
\[
\lim_{\alpha\to\infty} \| f(x) - J_\alpha f (x)     \|_{(PW,l^1)}=0
\]
Since
\[
\| f (x)  \|_{L^2(\R)} = \| \hat{f}(\xi)  \|_{L^2{\R}} \leq \sum_{j\in\Z}\|\hat{f}(\xi+2\pi j)  \|_{L^2(\T)}=\| f(x) \|_{(PW,l^1)},
\]
we also have that
\[
\lim_{\alpha\to\infty} \| f(x) - J_\alpha f (x) \|_{L^2(\R)} = 0.
\]
We move on now to the pointwise estimate.  Using the inversion formula and the Schwarz inequality as in the proof of theorem 2 in \cite{me}, we have
\begin{align*}
&\left| f(x)-J_\alpha f (x)   \right|  \\
\leq & C\sum_{j\in\Z}\| \hat{f}(\xi+2\pi j) - \sum_{m\in\Z}\widehat{I_\alpha g_m}(\xi-2 \pi m + 2\pi j)    \|_{L^2(\T)}
\end{align*}
We recognize this as \eqref{eq2}, and by the same reasoning used above, this tends to 0.  Since the estimates do not depend on the $x$ value, the convergence is uniform.
We summarize these ideas in the following theorem.
\begin{thm}
Given a function $f\in (PW,l^1)$, a complete interpolating sequence $\{x_n:n\in\Z\}$, and a regular family of interpolators $\{\phi_\alpha: \alpha\in A\}$, we have
\begin{enumerate}
\item[A.]  $\displaystyle \lim_{\alpha\to\infty}\| f- J_\alpha f  \|_{L^2(\R)}=0$,
\item[B.]  $\displaystyle \lim_{\alpha\to\infty}\| f- J_\alpha f  \|_{(PW,l^1)}=0$, and
\item[C.]  $\displaystyle \lim_{\alpha\to\infty}| f(x) - J_\alpha f (x)  | =0$, uniformly on $\R$.
\end{enumerate}
\end{thm}

\section{Examples and Conclusions}
Work of this type has been carried out before on the Paley-Wiener space $PW_\pi$.  The sequence of papers \cite{paper 3}, \cite{siva}, and \cite{me} all share this theme.  We list some specific examples that may be used to recover functions in $(PW,l^1)$.
\begin{enumerate}
\item the family of Gaussian kernels $\{ \exp(-x^2/4\alpha) : \alpha \geq 1\}$ 
\item the family of Poisson kernels $\{ \alpha/(x^2+\alpha^2) : \alpha \geq 1  \}$
\item the family of tempered splines whose degree increases to infinity
\end{enumerate}

The last item is not a regular family of interpolators, although it exhibits the same interpolation and recovery properties as regular families of interpolators so it may be used as well.  Examining the recovery property for $(PW,l^1)$, we find a sequence of coefficients $\{ a_{m,n}(\alpha):m,n\in\Z, \alpha\in A \}$ which satisfies
\[
\lim_{\alpha\to\infty}\sum_{m,n\in\Z}a_{m,n}e^{2\pi i m x}\phi_{\alpha}(x-x_n)=f(x).
\]
This looks quite a bit like a Gabor system, although the are important differences.  The space of applicability is reduced from $L^2(\R)$ to $(PW,l^1)$, which loses a good deal of structure due to the absence of an inner product.  However, this is a fairly large subspace, in fact, the following subspaces are contained in $(PW,l^1)$:
\begin{enumerate}
\item[a] the Schwartz space $\mathscr{S}$,
\item[b] $\displaystyle\bigcup_{\epsilon>0} \{ f\in L^2(\R) : \text{supp}(\hat{f})\subseteq [-\epsilon,\epsilon]    \}$, and
\item[c] $\displaystyle\bigcup_{\epsilon>0} \{ f\in L^2(\R) :  \hat{f}(\xi)= O (|\xi|^{-1-\epsilon})    \}$.
\end{enumerate}

This shows that the methods here are really an improvement of the results found in \cite{me}, \cite{paper 3}, and \cite{siva}, since this space is more than just a Paley-Wiener space.  It seems reasonable that these types of arguments can be extended out of the Hilbert space setting by using frame theory, but the calculations remain in progress.


\begin{thebibliography}{[00]}


\bibitem{Grochenig}
K. Gr{\"o}chenig.
\emph{Foundations of Time Frequency Analysis},
Birkh{\"a}user, Boston, MA, 2001. 

\bibitem{me}
J. Ledford.
``Recovery of Paley-Wiener functions using scattered translates of regular interpolators,"
J. Approx. Theory {\bf 173} (2013), 1-13.

\bibitem{paper 3}
Y. Lyubarskii \and W. Madych.
``Recovery of Irregularly Sampled Band Limited Functions via Tempered Splines," 
J.  Funct. Anal. \text{\bf{125}} (1994), 201-222.

\bibitem{siva}
T. Schlumprecht \and N. Sivakumar.
``On the Sampling and Recovery of Bandlimited Functions via Scattered Translates of the Gaussian,"
J. Approx. Theory \text{\bf{159}} (2009), 128-153.



\end{thebibliography}
\end{document}